\documentclass[a4paper,11pt]{amsart}
\usepackage{amsmath,inputenc,euscript,amssymb,geometry}
\usepackage{graphicx}
\usepackage{amssymb}
\usepackage{latexsym}
\usepackage{amssymb,amsbsy,amsmath,amsfonts,amssymb,amscd}
\usepackage{hyperref}
\usepackage{graphics}
\usepackage{color}
\newtheorem{lemma}{Lemma}[section]
\newtheorem{theorem}[lemma]{Theorem}

\begin{document}

\title 
[Evolution by VFE of curves with a corner]
{Evolution by the vortex filament equation\\ of curves with a corner}
\author[V. Banica]{Valeria Banica}

\address[V. Banica]{Laboratoire Analyse et probabilit\'es (EA 2172)\\D\'eptartement de Math\'ematiques\\ Universit\'e d'Evry, 23 Bd. de France, 91037 Evry\\ France, Valeria.Banica@univ-evry.fr}

\thanks{The author is partially supported by the French ANR project SchEq ANR-12-JS-0005-01}
\email{Valeria.Banica@univ-evry.fr}



\begin{abstract}
In this proceedings article we shall survey a series of results on the stability of self-similar solutions of the vortex filament equation. This equation is a geometric flow for curves in $\mathbb R^3$  and it is used as a model for the evolution of a vortex filament in fluid mechanics. The main theorem gives, under suitable assumptions, the existence and description of solutions generated by curves with a corner, for positive and negative times. Its companion theorem describes the evolution of perturbations of self-similar solutions up to a singularity formation in finite time, and beyond this time. We shall give a sketch of the proof. These results were obtained in collaboration with Luis Vega.

\end{abstract}

\maketitle

In the first section we shall present the vortex filament equation. Then in the second section we shall give the topics we are interested in, and state our results. The sketch of the proofs will be given in the last section.

\section{The Vortex Filament Equation}\label{sect:VFE}
The vortex filament equation (VFE) is the geometric flow of curves $\chi(t)$ in $\mathbb R^3$ governed by
\begin{equation}\label{VFE}\chi_t=\chi_x\wedge\chi_{xx}.
\end{equation}
Here $x$ stands for the arclength parameter of the curve $\chi(t)$. If one uses the Frenets system that characterize the tangent, normal and binormal vectors of a $\mathbb R^3$ curve in terms of the curvature and torsion, 
\begin{equation}\label{Frenet}
\left(\begin{array}{c}
T\\n\\b
\end{array}\right)_x=
\left(\begin{array}{ccc}
0 & c & 0 \\ -c & 0 & \tau \\ 0 &  -\tau & 0 
\end{array}\right)
\left(\begin{array}{c}
T\\n\\b
\end{array}\right),
\end{equation}
one gets that the equation (VFE) can be written as
\begin{equation}\label{binormal}\chi_t=c\,b.
\end{equation}
This emphasize that the curve evolve in the direction of the binormal vector, with speed proportional to the curvature. The equation is also known as the binormal flow. It was derived by Da Rios in 1906 and rediscovered by Arms and Hama in 1965 as a model for the dynamics in a 3D inhomogeneous inviscid Newtonian fluid of a vortex filament located at time $t$ on a $\mathbb R^3$ curve $\chi(t,x)$ (\cite{DaR,ArHa}, see also \cite{MaBe,Ri1,Ri2}). The approximation uses a truncation of Biot-Savart's law, which is also known as the Local Induction Approximation (LIA). 
The validity of this model has been proved for vorticities supported in an annulus (vortex rings), in the case of axial symmetry without swirl  (\cite{Ben-Cag-Mar}, see also \cite{Mar-rings, Mar-Neg} for related results). We recall here that the binormal flow is also conjectured to model the evolution of quantized vortex  
filaments in a Bose condensate in the incompressible limit. This was proved in the case of vortices concentrated around a circle, and a similar statement holds for higher dimensions (\cite{Je}). In a related setting, it was proved in \cite{BeOrSm} that mean curvature flow governs the dynamics in parabolic Ginzburg-Landau equation.

Although in general part of the complexity of the fluid equations might be lost through this approximation, this model is a simple and very rich one. For instance, it is a completely integrable equation, and this feature has been exploited to understand the topological properties of knotted vortex filaments (see \cite{Cal,Laf,Mag} and the references therein). It is worth pointing out also that the $\mathbb S^2$ tangent vector $T(t,x)$ satisfies the {Schr\"odinger map arising in ferromagnetism theory. Moreover, a link with the 1D cubic nonlinear Schr\"odinger equation (NLSE) is made at the second order of derivative in space as follows.

If one considers the filament function 
$$\psi(t,x)=c(t,x)e^{i\int_0^x\tau(t,s)ds},$$
it is easy to check that it satisfies the (NLSE) 
 \begin{equation}\label{Hasimoto}
i\psi_t+\psi_{xx}+\left(|\psi|^2-A(t)\right)\psi=0,
\end{equation}
wihere $A(t)$ is in terms of the curvature and torsion $(c,\tau)(t,0)$.
This fundamental remark has been made by Hasimoto in 1972 and it has allowed the transfer of informations from (NLSE) to (VFE) (\cite{Ha}). This Hasimoto transform can be seen as an inverse of the Madelung transform that connects the Gross-Pitaevskii equation to Euler equation with quantum pressure. Actually also here, the system satisfied by the curvature and torsion is a Euler-Korteweg one. 

The fact that the curvature is not allowed to vanish in order to define the filament function is just a technical obstruction. This was proved by Koiso in \cite{Ko} by considering another frame $(T,e_1,e_2)$ than the Frenet frame $(T,n,b)$, governed by
 \begin{equation}\label{Koiso}
\left(\begin{array}{c}
T\\e_1\\e_2
\end{array}\right)_x=
\left(\begin{array}{ccc}
0 & \alpha & \beta \\ -\alpha & 0 & 0\\ -\beta &  0 & 0 
\end{array}\right)
\left(\begin{array}{c}
T\\e_1\\e_2
\end{array}\right).
\end{equation}
If we denote $N=e_1+ie_2$ and if $T(t)$ is the tangent vector of a curve $\chi(t)$ solution of \eqref{VFE}, then one can compute
 \begin{equation}\label{derKoiso}
T_x=\Re(\overline \psi N), N_x=-\psi T, T_t=\Im\overline{\psi_x}N,\quad N_t=-i\psi_x T+i(|\psi|^2-A(t))\, N,
\end{equation}
and verify that
$$\psi(t,x)=\alpha(t,x)+i\beta(t,x)$$
is a solution to \eqref{Hasimoto} with $A(t)=\alpha^2(t,0)+\beta^2(t,0)$. Such a frame can be obtained by a rotation of the Frenet frame, $N(t,x)=(n+ib)(t,x)e^{i\int_0^x\tau(t,s)ds}$.

Conversely, given $\psi$ a solution of \eqref{Hasimoto}, $(e_0, e_1, e_2)$ an orthonormal basis of $\mathbb R^3$ and $P$ a point of $\mathbb R^3$, one might construct a solution of \eqref{VFE} in the following way. First, $(T,N)(t,x)$ can be constructed by imposing  $(T,N)(t_0,x_0)=(e_0,e_1+ie_2)$ and the evolutions laws \eqref{derKoiso}. Then $\chi(t,x)$ defined as 
\begin{equation}\label{chi}
\chi(t,x)=P+\int_t^{t_0}(T\wedge T_{xx})(\tau,x_0)d\tau+\int_x^{x_0}T(t,s)ds,
\end{equation}
is a solution of \eqref{VFE}. 
Summarizing, this recipe can be used to construct solutions of (VFE) starting from solutions of (NLSE). However, recovering the geometric properties of the solution of (VFE) is not obvious at all.  
It has been possible doing so starting from explicit solutions of (NLSE) as 
\begin{itemize}
  \item $\,\psi(t,x)=0,\, A(t)=0$.
  \item $\,\psi(t,x)=1,\, A(t)=-1$.
  \item $\psi(t,x)=e^{-itN^2} e^{iNx}, \, A(t)=-1.$
   \item$\psi(t,x)=e^{-itN^2} e^{iNx}\frac1{2\sqrt2}\, \frac1{\cosh(x-2Nt)},\,A(t)=-1.$
\end{itemize}
These cases give solutions to (VFE) that are still lines, circles traveling in the binormal direction, evolving helices and traveling waves respectively. It is important to note that the first three families of solutions are known dynamics of vortex filaments in fluids. However, the forth one was not known until Hasimoto's work in 1972 on (VFE). This was the starting point for the physicists Hopfinger and Browand that succeed to display such dynamics in a fluids experiment (\cite{HoBr}).

Well-posedness results for (VFE) for regular curves have been obtained by various methods \cite{DaR,LC,Ha,NiTa,TN}. One way is by using the link with (NLSE), having in mind that this equation is well posed in $H^s$, $s\geq 0$ on the line and on the circle \cite{CW,Bo}. Recently the less regular case of currents has been considered in a weak formulation of the equation \cite{JeSm1,JeSm2}.\\

We shall focus now on the self similar solutions of (VFE), that is solutions of the type 
$$\chi(t,x)=\sqrt{t}G\left(\frac x{\sqrt{t}}\right).$$
We recall that since the 70's the (VFE) and its self-similar solutions were considered in works on vortex dynamics in superfluids \cite{Sc,Bu,Li}, in ferromagnetism \cite{LRT,LD}, in aortic heart valve leaflet miocardic modeling \cite{PMQ,V},.

In \cite{LRT} it was shown that self-similar solutions form a family $\{\chi_a\}_{a\in\mathbb R^{+*}}$ caracterized by the explicit curvature and torsion 
$(c_a,\tau_a)(t,x)=\left(\frac{a}{\sqrt{t}},\frac{x}{2t}\right)$. Numerical computations on this formation of a singularity in finite time were given in \cite{Bu}. In \cite{GRV} this was proved rigorously. More precisely, it was shown that a corner appear at time $t=0$:
$$\left|\chi_a(t,x)-x(A^+_a  \mathbb{I}_{[0,\infty)}(x)+A^-_a  \mathbb{I}_{(-\infty,0]}(x))\right|\leq {2a}{\sqrt{t}},$$
with $A^\pm_a\in\mathbb S^2$ distincts, non-opposite and $\sin\frac{(\widehat{A_a^+,-A_a^-})}{2}=e^{-\frac{a^2}{2}}.$ In particular, any corner can be obtained in finite time from a rotated and translated self-similar solution. 

At the level of vortex filament dynamics in fluids, an analogy can be done between the evolution of $\chi_a$ and the ``delta wing" vortex. More precisely, in \cite{HGCV} numerical simulations for selfsimilar solutions of the binormal flow are in correlation with the physical experiment. 
\begin{center}
\includegraphics[width=3.6in]{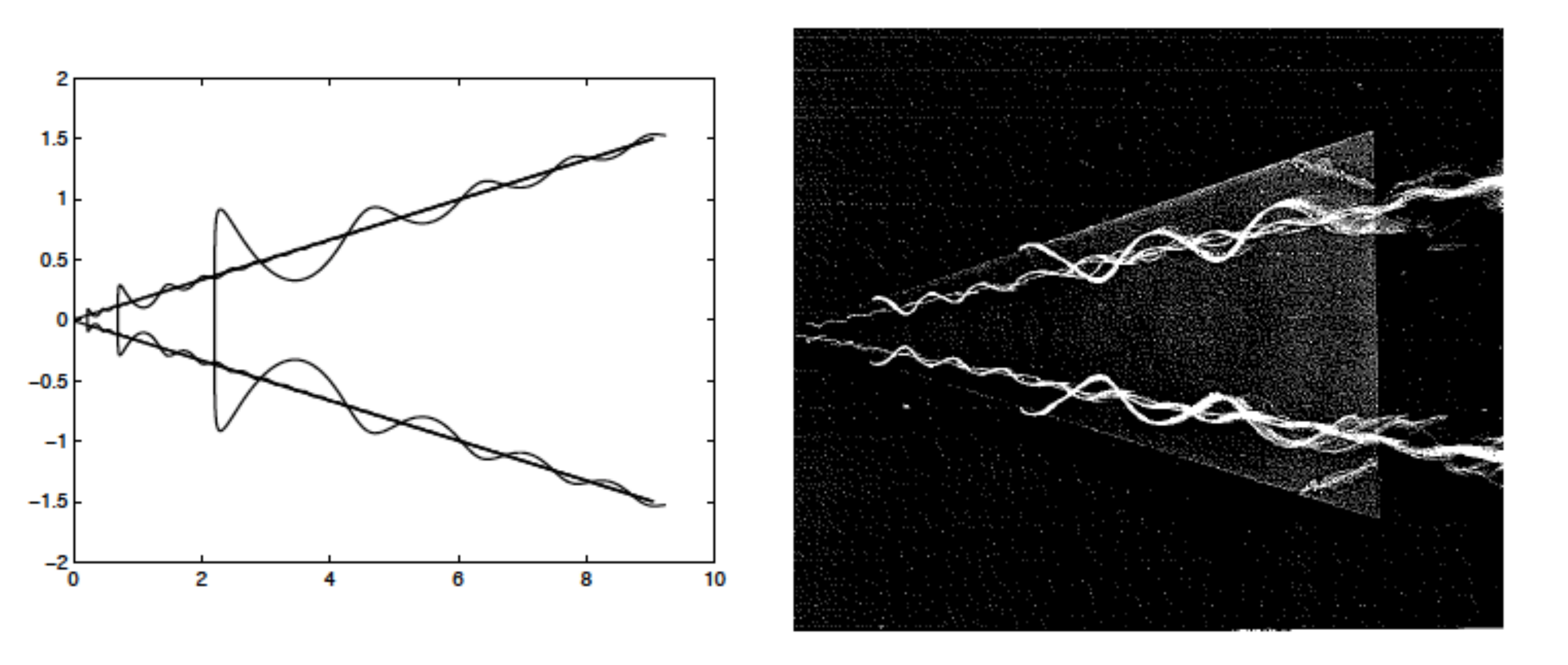}\\ On the left the numerical simulations of $\chi_a$ from \cite{HGCV}; on the right the ``delta wing" vortex formed by fluid flowing over a triangular obstacle.
\end{center}

\section{Statements of the results}

The question that motivated our series of papers \cite{BV1,BV2,BV3,BV4} is whether the formation of singularity in finite time for self-similar solutions is stable or not. For instance, if one consider a small perturbation of $\chi_a(1)$ and let it evolve through (VFE), does it generate also a singularity in finite time, and if it is the case, which is the geometrical description of the singularity. In \cite{BV1,BV2,BV3} we have answered these questions and this allowed us to treat in \cite{BV4} the issue of considering as initial data for (VFE) curves having a corner.

The first obstacle in studying perturbations of $\chi_a$ is that this particular solution does not enter the framework of local well-posedness results for the binormal flow. Moreover, it is related to the delicate issue of rough data for the 1D cubic Schr\"odinger equation, as follows. Via Hasimoto's transform the behaviour at $t=0$ of perturbations of $\chi_a$, solutions of the binormal flow, can be understood from the behaviour at $t=0$ of perturbations of its filament function 
$$\psi_a(t,x)=\frac{a}{\sqrt{t}}\,e^{i\frac{x^2}{4t}},$$
solutions of 
 \begin{equation}\label{NLS}i\psi_t+\psi_{xx}+\frac{1}{2}\left(|\psi|^2-\frac{a^2}{t}\right)\psi=0.\end{equation}
The corner $\chi_a(0,x)$ corresponds to $\psi_a(0,x)=a\,\delta_{x=0}$. 

This calls for a review on the Cauchy problem for the 1D cubic Schr\"odinger equation: 
$$
 i\psi_t+\psi_{xx}\pm |\psi|^2\psi=0.
$$
 This equation is well-posed in $H^{s}$, $s\geq 0$ \cite{GV,CW}. If $s<0$ it is ill-posed \cite{KPV,CCT}. Well-posedness holds for data with summable Fourier transform \cite{VV,Gr}. For initial data precisely $a\delta_{x=0}$, the equation is ill-posed: if there is uniqueness, the solution is $e^{ia^2\log t}\psi_a(t,x)$, which does not have a limit at $t=0$ \cite{KPV}. A natural change to avoid the logarithmic phase is to consider the perturbed equation \eqref{NLS} together with its solution $\psi_a(t,x)$.  
Note that for the free equation, $\psi_a(t,x)$ is a solution with initial data $a\delta_{x=0}$, and smooth perturbations of it at time $t=1$ evolve to smooth perturbations of $a\delta_{x=0}$ at time $t=0$. In \cite{BV2} we have proved that small smooth perturbations of the solution $\psi_a(t,x)$ at time $t=1$, evolving through \eqref{NLS} still do not have a limit at $t=0$.
This should be seen as a "blow-up" result for 1D cubic Schr\"odinger equation with rough data.\\

In order to understand the behavior at time $t=0$ of perturbations $\psi$ of $\psi_a$, solutions of \eqref{NLS}, we shall use the pseudo-conformal transformation. More precisely, defining 
$$\psi(t,x)=\frac{1}{\sqrt{t}}\,e^{i\frac{x^2}{4t}}\,\overline{v}\left(\frac 1t,\frac xt\right)$$ 
the problem reduces to large time behavior of perturbations of $v_a(t,x)=a$, solutions of:
$$iv_t+v_{xx}+\frac 1{2t}\left(|v|^2-a^2\right)v=0.$$
Setting $u=v-a$ we are now interested in the large time behavior of small data for the equation
\begin{equation}\label{NLSu}iu_t+u_{xx}+\frac 1{2t}\left(|u+a|^2-a^2\right)(u+a)=0.\end{equation}
Here we note that the zero Fourier mode of $u(t)$ turns out to grow logarithmically in time for generic data (see Appendix B of \cite{BV2}).

We shall prove scattering results for this problem, in the sense of linking the behavior at large times of nonlinear solutions to the one of linear evolutions. As usual, the worse terms to control  when trying to prove scattering results are the lowest order terms. In the case of \eqref{NLSu} these are two linear terms: $\frac{a^2}{2t}u+\frac{a^2}{2t}\,\overline{u}$. The potential $\frac{a^2}{2t}$ is not integrable, so we shall discard the non-oscillant term $\frac{a^2}{2t}u$.  For small data, we shall get a long-range asymptotic profile for $u(t)$: 
\begin{equation}\label{longrange}e^{i a^2\log \sqrt{t}}\,e^{it\partial_x^2}\,u_+(x).\end{equation}
Having in mind that linear evolutions in $\mathbb R^n$ decay in time like $t^{-\frac n2}$, the present situation could be linked to the scattering results obtained for the 1D cubic Schr\"odinger equation \cite{Oz,Ca,HN}, the 2D Gross-Pitaevskii equation \cite{GNT}, or the 2D quadratic Schr\"odinger equation \cite{MTT,ST,GMS}. However, the proof will require a specific analysis of the linear equation and the introduction of new appropriate spaces. 

Let us point out that if $xu_+\in L^2$, by getting back through the pseudo-conformal transform to the filament function, we obtain that it behaves like: 
$$\frac{a}{\sqrt{t}}\,e^{i\frac{x^2}{4t}}+\frac{e^{ia^2\log \sqrt{t}}}{\sqrt{4\pi i}}\widehat{\overline{u_+}}\left(-\frac x2\right).
$$
This shows that indeed small smooth perturbations of $\psi_a$ at time $t=1$, evolving through \eqref{NLS} still do not have a limit at $t=0$. Nevertheless, this will not be an obstruction to have a limit at $t=0$ for the corresponding perturbations of $\chi_a$, solutions of (VFE). This is not contradictory since the filament function is defined in terms of quantities appearing at the second order derivative in space of the curve. \\ 

Now we shall review the previous results on singularity formation obtained in \cite{BV1,BV2,BV3}. Each article contains a result at the (NLSE) level and informations on the formation of singularities at the (VFE). The passage from (NLSE) to (VFE) in \cite{BV1,BV3} is particularly lengthy. We summarize these results as follows:
\begin{itemize}
\item (\cite{BV1}): if $s\in\mathbb N$ there exist modified wave operators in $H^s$ for $u_+$ small in $\dot{H}^{-2}\cap H^s\cap W^{s,1}$ (i.e. there exists a nonlinear solution $u$ of \eqref{NLSu} that behaves at large times in $H^s$ like \eqref{longrange}).\\ If moreover $u_+$ is small in weighted spaces, and $s\geq 3, $ there exists solutions $\chi$ of (VFE) with corner-type singularity (i.e. the formation of singularities for self-similar solution is not an isolated phenomena).

\item (\cite{BV2}): if $s\in\mathbb N$ asymptotic completeness holds (with weaker decay than in \cite{BV1}) in $H^s$ for $\partial_x^k u(1)$ small in 
$$X^\gamma=\{f\in L^2, |\xi|^\gamma \hat f(\xi)\in  L^\infty(|\xi|\leq 1)\},$$
$0\leq k\leq s,\gamma<\frac 12$, with final state $u_+\in X^{\gamma^+}$ (i.e. all initial data $u(1)$ will have a global evolution through \eqref{NLSu} and will behave like \eqref{longrange} for some $u_+\in X^{\gamma^+}$). Moreover there exists modified wave operators for final data in $ X^{\gamma^+}$, so the scattering operator can be defined.\\
If $s\geq 3$ all small perturbations of $\chi_a$ at time $t=1$ generate a singularity at time $t=0$ (i.e. stronger stability of the singularity formation as in \cite{BV1}, but with weaker geometric description of the singularity formation).

\item (\cite{BV3}): If moreover $u(1)$ small in weighted spaces, $u(t)$ is in weighted spaces with a good control in time and we get strong stability results on the singularity (same corner for instance).
 \end{itemize}

In view of these results, the natural question is to understand the extension after singularity time of the perturbations of $\chi_a$ at time $t=1$.  
Unfortunately the weighted spaces do not fit for a scattering functional setting at the (NLSE) level. This issue has been avoided in \cite{BV4}, by using the scattering operator in $X^{\gamma^+}$ from \cite{BV2} and by removing the weights hypothesis in \cite{BV3}. One of the two companion results in \cite{BV4} is the  following.

\begin{theorem}\label{the}(Continuation of solutions after singularity time, \cite{BV4})
Let $\chi(1)$ be a small perturbation of a self-similar solution $\chi_a$ at time $t=1$ in the sense that the filament function of $\chi(1)$ is $(a+u(1,x))e^{i\frac{x^2}{4}}$, with $\partial_x^ku(1)$ small in $X^\gamma$ with respect to $a$ for all $0\leq k\leq 4$, for some $\gamma<\frac 12$. \\
\medskip
We construct a solution $\chi\in\mathcal C([-1,1],Lip)\cap \mathcal C([-1,1]\backslash \{0\},\mathcal C^4)$ for the binormal flow on $ t\in [-1,1]\backslash \{0\}$, which is a weak solution on the whole interval $[-1,1]$. The solution $\chi$ is unique in the subset of $\mathcal C([-1,1],Lip)\cap \mathcal C([-1,1]\backslash \{0\},\mathcal C^4)$ such that the associated filament functions at times $\pm1$ can be written as $(a+u(\pm1,x))e^{i\frac{x^2}{4}}$ with $\partial_x^ku(\pm1)$ small in $X^\gamma$ with respect to $a$ for all $0\leq k\leq 4$. 

\medskip
\noindent
Moreover,the solution $\chi$ enjoys the following properties:
\begin{itemize}

\item there exists a limit of $\chi(t,x)$ and of its tangent vector $T(t,x)$ at time zero, and
$$\sup_s|\chi(t,x)-\chi(0,x)|\leq C\sqrt{|t|},\quad\sup_{|x|\geq\epsilon>0}|T(t,x)-T(0,x)|\leq C_\epsilon |t|^{\frac 16^-},$$

\item $\forall t_1, t_2\in[-1,1]\backslash \{0\}$ the following asymptotic properties hold  
$$ \chi(t_1,x)-\chi(t_2,x)=\mathcal O\left(\frac{1}{x}\right),\qquad T(t_1,x)- T(t_2,x)=\mathcal O\left(\frac{1}{x}\right),$$

\item $\exists T^\infty\in\mathbb S^2, N^\infty\in\mathbb C^3$ such that uniformly in $-1\leq t\leq 1$,
$$T(t,x)-T^\infty=\mathcal O\left(\frac{1}{\sqrt x}\right),\quad (n+ib)(t,x)-N^\infty\,e^{ia^2\log\frac{\sqrt t}{x}-i\frac{x^2}{4t}}=\mathcal O\left(\frac{1}{\sqrt x}\right),$$
\item modulo a rotation and a translation, we recover at the singularity point $(0,0)$ the same structure as for $\chi_a$:
$$\lim_{x\rightarrow 0^\pm}T(0,x)=A_a^\pm\quad,\quad  \lim_{x\rightarrow 0^\pm} \lim_{t\rightarrow 0} (n+ib)(t,x)e^{-ia^2\log\frac{\sqrt t}{x}+i\frac{x^2}{4t}}=B_a^\pm.$$
\end{itemize}
\end{theorem}

We insert here the picture of continuation through time $t=0$ in the very particular symmetric case of self-similar solutions $\chi_a$:

\begin{center}
\includegraphics[width=1.5in]{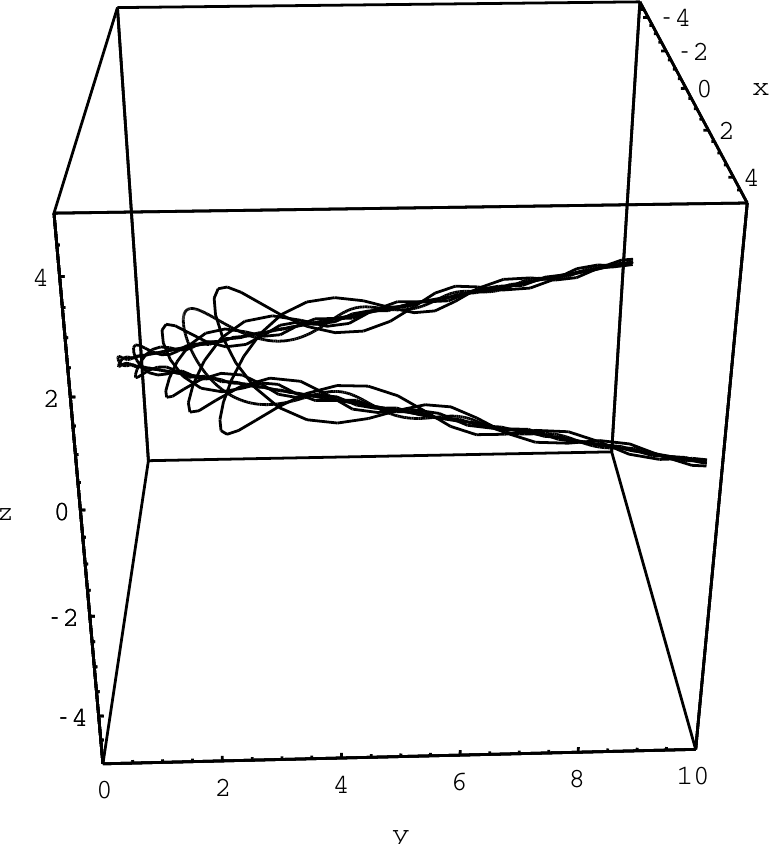}
\end{center}
This is obtained simply by using the time-reversibility of (VFE). It yields that the extension on negative times must be a solution for positive time with initial data $\chi_a(0,-x)$. Such a solution is given precisely by a $\pi-$rotation of $\chi_a$ with respect to the bisector of the corner. In the general case considered in the theorem, the extension through time $t=0$ is very involved.

We end with the main theorem in \cite{BV4} (its proof is a main step for Theorem \ref{the}, that we have preferred to state first for presentation reasons). 
\begin{theorem}\label{theivp}(The binormal flow with initial values with a corner, \cite{BV4})
Let $\chi_0$ be a smooth $\mathcal C^4$ curve, except at $\chi_0(0)=0$ where a corner is located, i.e. that there exist $A^+$ and $A^-$ two distinct non-colinear unitary vectors in $\mathbb R^3$ such that
$$\chi_0'(0^+)=A^+,\quad \chi_0'(0^-)=A^-.$$
We set $a$ to be the real number given by the unique self-similar solution of the binormal flow with the same corner as $\chi_0$ at time $t=0$. 
We suppose that the curvature of $\chi_0(x)$ (for $x\neq 0$)  satisfies $(1+|x|^4)c(x)\in L^2$ and $|x|^{\gamma}c(x)\in L^\infty_{(|x|\leq 1)}$, small with respect to $a$. 
%We denote $T_0=\chi_0'$. We define for $x>0$  a complex-valued function $g$ and a $\mathbb C^3$-valued function $\tilde N_0$ by solving  the system  $$ T_{0x}(x)=-\Re (g(x)\tilde N_0(x)),\,\tilde N_{0x}(x)=-\overline{g}(x)T_0(x),\, (T_0,\tilde N_0)(0)=(A^+,B^+) $$ We define $g(x)$ and $\tilde N_0$ similarly for $x<0$ imposing $(A^-,B^-)$ as initial data.  We suppose $\chi_0$ to be such that $(1+|x|^4)g\in L^2$ and $\frac{g(x)}{|x|^{2\gamma}}\in L^\infty_{(x^2\leq 1)}$ are small with respect to $a$ for some $0<\gamma<\frac 14$. 
\medskip

\noindent
Then there exists 
$\chi(t,x)\in\mathcal C([-1,1],Lip)\cap \mathcal C([-1,1]\backslash \{0\},\mathcal C^4),$
regular solution of the binormal flow for $t\in [-1,1]\backslash \{0\}$, having $\chi_0$ as value at time $t=0$, and enjoying all properties from the previous theorem.
\end{theorem}

The ordered main steps of the proof of Theorem \ref{the} (including the proof of Theorem \ref{theivp}) are the following:
\begin{itemize}
\item Asymptotic completeness for the solution $u$ of \eqref{NLSu} with initial data $u(1)$ small in $X^\gamma$.
\item Asymptotics in space and in time for the tangent and normal vectors of the curve $\chi$ constructed from $u$ as explained in the first section.
\item Existence of a trace of the tangent and modulated normal vectors at time $t=0$.
\item Recovering of the self-similar singularity structure at $t=0,x=0$ for $\chi$ (for instance same angle as $\chi_a$).
\item Recipe for the IVP for curves with data $\chi_0$ with a corner (involves the existence of modified wave operator for the same equation \eqref{NLSu} with final data $g$, a function defined from a system of the traces of tangent and normal vectors of $\chi_0$). This ends the proof of the main part of Theorem \ref{theivp}.
\item Construction of $\chi$ for negative times. 
\end{itemize}

\section{Sketch of the proof}
In the first subsection we shall present the scattering results for \eqref{NLSu}. From a long-range solution of \eqref{NLSu} we construct a solution $\chi$ of \eqref{VFE} as explained in the first section. In \S \ref{subsect:tang} we shall derive asymptotics for its tangent and normal vectors in time and space and we shall obtain their traces at time $t=0$. Finally, in the last subsection we shall start by describing the structure of the singularity, which will allow to obtain the recipe of the IVP, and then to extend the solution $\chi$ through time $t=0$.

\subsection{Scattering for cubic 1DNLS} 
As mentioned in the previous section, we make a phase change in \eqref{NLSu} to get rid of the non-oscillant term $\frac{a^2}{2t}u$. By changing $u(t)$ into $u(t)e^{- ia^2\log t}$ the new equation writes 
\begin{equation}\label{NLSuphased} iu_t+u_{xx}\pm\frac{a^2}{t^{1\pm 2ia^2}}\,\overline{u}=\frac{u^{2,3}}{t}.\end{equation}
We shall first treat the linear equation 
\begin{equation}\label{linNLSuphased}iu_t+u_{xx}\pm\frac{a^2}{t^{1\pm 2ia^2}}\,\overline{u}
=0,\end{equation}
and eventually deal with the nonlinear one by a perturbative argument.

We start with a paragraph on a-priori bounds at the linear level and by pointing out that there is a growth of the low Fourier modes. This will give us a hint for finding the appropriate space for global existence and scattering, that we shall present in the following paragraph.

\subsubsection{Growth of the Fourier modes for the linear equation}
For the linear equation \eqref{linNLSuphased} we have 
$$\partial_t\frac{\|u(t)\|_{L^2}^2}{2}\leq \frac{a^2}{t}\|u(t)\|_{L^2},$$
so we get a polynomial control 
 $$\|u(t)\|_{L^2}\leq t^{a^2}\|u(1)\|_{L^2}.$$
 In particular we obtain the existence of an unique solution $u\in\mathcal C([1,\infty[,L^2).$ For our scattering purposes we would like to have $u\in L^\infty ([1,\infty[,L^2).$ So we start an analysis in Fourier variables. 
As a first a-priori bound we get, as for the $L^2$ norm,
$$|\hat{u}(t,\xi)|\leq \, t^{a^2} \,\left(|\hat{u}(1,\xi)|+|\hat{u}(1,-\xi)|\right).$$
However, by a precise lengthy study of the equations on $\widehat{\Re u},\widehat{\Im u}$ we improve to  (see \S 2.1 in \cite{BV2}, and \S 6.1 in \cite{BV3}):
\begin{equation}\label{est1}|\hat{u}(t,\xi)|\leq \, C(\delta)\,t^{\delta} \,\left(|\hat{u}(1,\xi)|+|\hat{u}(1,-\xi)|\right),\end{equation}
\begin{equation}\label{est2}|\hat{u}(t,\xi)|\leq \, \left(1+\frac{C(\delta)}{\xi^{\delta}}\right) \,\left(|\hat{u}(1,\xi)|+|\hat{u}(1,-\xi)|\right), \,\forall \delta>0.\end{equation}
In particular, this means that for $|\xi|\geq 1$ the Fourier mode has bounded variation in time, while for $|\xi|< 1$ there might be a growth. Actually time logarithmic growth for $\xi=0$ is observed for well-chosen data, and this also happens at nonlinear level for equation \eqref{NLSuphased} (see Appendix B of \cite{BV2}).

\subsubsection{Global existence and scattering}
In view of estimates \eqref{NLSuphased}-\eqref{linNLSuphased} we introduce for $\gamma,\tilde\gamma<\frac 12$ the spaces 
$$\|f\|_{X}=\|f\|_{L^2}+\||\xi|^{\gamma}\hat{f}(\xi)\|_{L^\infty(|\xi|\leq 1)},\qquad\qquad\qquad\qquad$$
$$\|g\|_{Y}=\sup_{t\geq 1}\,\left(\|g(t)\|_{L^2}+\frac{1}{t^{\tilde\gamma}}\||\xi|^{\gamma}\hat{g}(t,\xi)\|_{L^\infty(|\xi|\leq 1)}\right).$$
Global existence for the linear equation \eqref{linNLSuphased} in $Y$ follows easily for data in $X$ since on the one hand 
$$\frac{1}{t^{\tilde\gamma}}\||\xi|^{\gamma}\hat{u}(t,\xi)\|_{L^\infty(|\xi|\leq 1)}\leq C(\tilde\gamma)\||\xi|^{\gamma}\hat{u}(1,\xi)\|_{L^\infty(|\xi|\leq 1)}\leq  C(\tilde\gamma)\|u(1)\|_{X}.$$
and on the other hand
$$\|u(t)\|_{L^2}=\|\hat{u}(t)\|_{L^2(|\xi|\leq 1)}+\|\hat{u}(t)\|_{L^2(1\leq |\xi|)}$$
$$\leq C(\delta)\||\xi|^{-(\delta+\gamma)}\|_{L^2(|\xi|\leq 1)}\,\||\xi|^{\gamma}\hat{u}(1,\xi)\|_{L^\infty(|\xi|\leq 1)}+\|\hat{u}(1)\|_{L^2(1\leq |\xi|)}\leq  \tilde C(\delta)\|u(1)\|_{X}.$$
Therefore we have that $u$ is a global solution in $Y\subset L^\infty ([1,\infty[,L^2)$.\\

Asymptotic completeness for the linear equation \eqref{linNLSuphased} holds if and only if $e^{-it\Delta}u(t)$ has a limit in $L^2$ as $t$ goes to $\infty$. Therefore we need
$$\ \|e^{-it_1\Delta}u(t_1)-e^{-it_2\Delta}u(t_2)\|_{L^2}\overset{t_1,t_2\rightarrow\infty}{\longrightarrow}0,$$
which is equivalent, using the Duhamel formulation, to 
$$\left\|\int_{t_1}^{t_2}e^{-is\Delta}\overline{u}(s,x)\,\frac{ds}{s^{1\pm2ia^2}}\right\|_{L^2}\overset{t_1,t_2\rightarrow\infty}{\longrightarrow}0.$$
We estimate the Fourier modes of this Duhamel term by integrations by parts:
$$A_{t_1,t_2}(\xi)=\int_{t_1}^{t_2}e^{is\xi^2}\overline{\hat u(s,-\xi)}\,\frac{ds}{s^{1\pm2ia^2}}$$
$$=\left.\frac{e^{is\xi^2}}{i\xi^2}\frac{\overline{\hat u(s,-\xi)}}{s^{1\pm2ia^2}}\right|_{t_1}^{t_2}-\int_{t_1}^{t_2}\frac{e^{is\xi^2}}{i\xi^2}\left(\frac{\partial_s\overline{\hat u(s,-\xi)}}{s^{1\pm2ia^2}}-(1\pm 2ia^2)\frac{\overline{\hat u(s,-\xi)}}{s^{2\pm2ia^2}}\right)\,ds.$$
The term involving $\partial_s\overline{\hat u(s,-\xi)}$ has a non-integrable factor; we use the equation to express it as
$$\partial_s\overline{\hat u(s,-\xi)}=i\xi^2\overline{\hat u(s,-\xi)}\pm i\frac{a^2}{s^{1\mp2ia^2}}\hat u(s,\xi).$$
Notice that the second term will give an $1/s^2$ integrable factor, while the first one recovers $-A_{t_1,t_2}(\xi)$. This is due to the fact that the linear term is conjugated. We obtain 
$$2|A_{t_1,t_2}(\xi)|\leq \left.\frac{|\hat u(s,\xi)|}{s \,\xi^2}\right|_{t_1}^{t_2}+\int_{t_1}^{t_2}\frac{|\hat u(s,\xi)|+|\hat u(s,-\xi)|}{s^2\,\xi^2}\,ds.$$
Now we use \eqref{est2} to get
\begin{equation}\label{estA} 2|A_{t_1,t_2}(\xi)|\leq \left(1+\frac{C(\delta)}{\xi^{\delta}}\right) \,\frac{|\hat{u}(1,\xi)|+|\hat{u}(1,-\xi)|}{t_1\,\xi^2}, \,\,\forall \delta>0.\end{equation}
This immediately gives the estimate of the $L^2$ norm in two regions in Fourier variables :
$$\|A_{t_1,t_2}\|_{L^2(1\leq |\xi|)}\leq\frac{C}{t_1}\|u(1)\|_{L^2},$$
and
 $$\|A_{t_1,t_2}\|_{L^2(\frac{1}{t_1}\leq\xi^2\leq 1)}\leq\frac{\||\xi|^{\gamma}\hat u(1)\|_{L^\infty(|\xi|\leq 1)}}{t_1} \||\xi|^{-2-(\gamma+\delta)}\|_{L^2(\frac{1}{t_1}\leq|\xi|\leq 1)}\\\,\\
\leq \frac{C}{t_1^{\frac 14-\frac{\gamma+\delta}2}}\||\xi|^{\gamma}\hat u(1)\|_{L^\infty(|\xi|\leq 1)}.$$
On the region $ \frac{1}{t_2}\leq \xi^2\leq \frac{1}{t_2}$ we split the Duhamel integral into two parts and estimate them using \eqref{est2} and \eqref{estA} as follows:
$$|A_{t_1,t_2}(\xi)|\leq \int_{t_1}^{1/\xi^2}|\hat u(s,-\xi)|\,\frac{ds}{s}+|A_{1/\xi^2,t_2}(\xi)|\leq \frac{\|\hat u(1)\|_{L^\infty(|\xi|\leq 1)}}{|\xi|^{\delta}}|\log\xi^2|+\frac{\|\hat u(1)\|_{L^\infty(|\xi|\leq 1)}}{|\xi|^{\delta}}.$$
Therefore
$$\|A_{t_1,t_2}\|_{L^2(\frac{1}{t_2}\leq\xi^2\leq \frac{1}{t_1})}\leq \||\xi|^{\gamma}\hat u(1)\|_{L^\infty(|\xi|\leq 1)}\|\log\xi^2|\xi|^{-(\gamma+\delta)}\|_{L^2(\frac{1}{t_2}\leq\xi^2\leq \frac{1}{t_1})}\leq \frac{C}{t_1^{\frac 14-\frac{\gamma+\delta}2}}\||\xi|^{\gamma}\hat u(1)\|_{L^\infty(|\xi|\leq 1)}.$$
Finally, on $\xi^2\leq \frac{1}{t_2}$ we use \eqref{est2}
$$|A_{t_1,t_2}(\xi)|\leq \frac{\|\hat u(1)\|_{L^\infty(|\xi|\leq 1)}}{|\xi|^{\delta}}|\log t_2|\leq \frac{\|\hat u(1)\|_{L^\infty(|\xi|\leq 1)}}{|\xi|^{\delta}}|\log\xi^2|$$
and integrate
$$\|A_{t_1,t_2}\|_{L^2(\xi^2\leq \frac{1}{t_2})}\leq \frac{C}{t_2^{\frac 14-\frac{\gamma+\delta}2}}\||\xi|^{\gamma}\hat u(1)\|_{L^\infty(|\xi|\leq 1)}.$$
Summarizing, we have indeed 
$$\left\|A_{t_1,t_2}\right\|_{L^2}\overset{t_1,t_2\rightarrow\infty}{\longrightarrow}0,$$
so asymptotic completeness follows for \eqref{linNLSuphased}:
$$\exists \,u_+\in L^2,\, \|u(t)-e^{\pm i a^2\log t}\,e^{it\partial_x^2}u_+\|_{L^2}=\mathcal{O}(t^{-\frac 14+\frac{\gamma^+}2}).$$
Moreover, a posteriori estimates show that the linear solution $u$ belongs to $L^4((1,\infty),L^\infty)$. For proving this the worse term to estimate is
$$ia^2\int_1^te^{i(t-s)\Delta}\overline{e^{it\Delta} u_+}\frac{ds}{s^{1\pm 2ia^2}},$$
and it is treated by splitting it in the Fourier variable into the regions $\xi^2\leq\frac{1}{t}$, $\frac 1t\leq\xi^2\leq 1$ and $1\leq \xi^2$ (see \S2.4 in \cite{BV2} for details).

\bigskip
The nonlinear equation \eqref{NLSuphased} is treated as a perturbation of the linear one \eqref{linNLSuphased}, by a fixed point argument in $Y\cap L^4((1,\infty),L^\infty)$, and global existence and scattering in $L^2$ follow. Moreover, the final state belongs to $X^{\gamma^-}$.  
The analysis works the same for higher Sobolev spaces (see \S3 in \cite{BV2}).

\subsection{Existence and properties of $T(0,x)$}\label{subsect:tang} 
As stated in Theorem \ref{the}, we start with $u(1)$ small in $X$ together with its first four derivatives in space. From the previous subsection we obtain that there exists $u$ solution of \eqref{NLSu} that has the long-range asymptotic
$$e^{i a^2\log \sqrt{t}}\,e^{it\partial_x^2}\,u_+(x),$$
for some $u_+\in X^{\gamma^-}$. Now we define $\psi$ by the pseudo-conformal transform,
$$\psi(t,x)=\frac{e^{i\frac{x^2}{4t}}}{\sqrt{t}}\,(a+\overline{u})\left(\frac 1t,\frac xt\right),$$ 
and we define vectors $(T,N)(t,x)$ by imposing their derivatives to satisfy \eqref{derKoiso} and $(T,N)(1,0)$ to be the canonical basis of $\mathbb R^3$. Then we define $\chi$ by \eqref{chi} and $P=(0,0,0)$, for $t> 0$. In the next two paragraphs we shall prove that there is a trace at time $t=0$ for the tangent vector $T(t,x)$ of $\chi(t,x)$. The extension of $\chi$ for negative times, together with the initial value problem result will be reviewed in \S\ref{sect:IVP}.

\subsubsection{Asymptotics in space and in time for the tangent and normal vectors}
For fixed $t\in ]0,1]$ we look for the asymptotic in space of $T(t,x)$. Using \eqref{derKoiso} we have 
$$T(t,x_2)-T(t,x_1)=\Re \int_{x_1}^{x_2}\frac{e^{-i\frac{s^2}{4t}}}{\sqrt t}(a+ u)\left(\frac 1t,\frac st\right) N (t,s)\,ds.$$
We integrate by parts from the oscillatory phase, and by using again \eqref{derKoiso}, 
$$\left|T(t,x_2)-T(t,x_1)+\Im \int_{x_1}^{x_2} e^{-i\frac{s^2}{4t}}\frac{2\sqrt t}{s}u_s\left(\frac 1t,\frac st\right) N (t,s)ds\right|\leq C(u(1))\frac{\sqrt t}{x_1}.$$
Therefore we obtain the existence of a limit in space:
$$\exists T^\infty(t)=\underset{x\rightarrow\infty}{\lim}T(t,x).$$ 
Now, by integrations by parts in time using \eqref{derKoiso} (see \S 3 in \cite{BV3}) we get that the $T^\infty(t)$ is independent in time:
$$T^\infty(t)=T^\infty,\, \forall 0<t\leq 1,$$
so we have obtained
\begin{equation}\label{estT}\left|T(t,x)-T^\infty-\Im \int_{x}^{\infty} e^{-i\frac{s^2}{4t}}\frac{2\sqrt t}{s}u_s\left(\frac 1t,\frac st\right) N (t,s)ds\right|\leq C(u(1))\frac{\sqrt t}{x}.\end{equation}
In view of this formula, to understand better the tangent vector we need informations on the complex normal vector $N(t,x)$. The later is very oscillating, and we modulate it as follows:
$$\tilde N(t,x)=N(t,x)e^{i\Phi}\,\,\,,\,\,\, \Phi(t,x)=-a^2\log \sqrt{t}+a^2\log x.$$ 
Again by integrations by parts in space and in time (see \S 3 in \cite{BV3}), and a compactness argument (see \S 3.1 in \cite{BV4}) we obtain the existence of $N^\infty\in\mathbb C$ such that 
\begin{equation}\label{estN} \left|\tilde N(t,x)-N^\infty-i\int_x^\infty \overline{h}(t,s)T(t,s)\,ds\right|\leq C(u(1))\left(\frac{\sqrt t}{x}+\frac{t}{x^2}+\sqrt t\right).
\end{equation}
with
$$h(t,s)=\frac{2\sqrt t}{s}u_s\left(\frac 1t,\frac st\right) e^{-i\frac{s^2}{4t}}e^{-i\Phi}.$$
Moreover, starting from \eqref{estT} and \eqref{estN} we obtain self-similar estimate involving only $T$, $T^\infty$, $N^\infty$ and $u$:
$$T(t,x)-T^\infty -\Im N^\infty \int_x^\infty h(t,s)ds-\Re\int_x^\infty h(t,s)\int_s^\infty \overline{h(t,s')}T(t,s')ds'ds,$$
and this process is iterated by handling the multiple integrals of $h$ to get again self-similar estimates:
$$|T(t,x)-\sum_{j=1}^{2n}a_j(t,x)|\leq C(u(1))\left(\sqrt{t}+\frac{\sqrt{t}}{x}+\frac{t\sqrt{t}}{x^3}\right),$$
with $a_j(t,x)$ explicit multiple integrals involving $h$ (see Lemma 3.3 in \cite{BV4}).

\subsubsection{The limit in time for $T(t,x)$ at fixed $x\neq 0$}
By denoting $\tilde h(s)=\frac{i\widehat{u_+}\left(\frac s2\right) }{s^{ia^2}}$ we prove that for $0<t\leq 1$ and $0<x$ (see Lemma 3.5 in \cite{BV4}):
$$
\left|\int_x^\infty h(t,s_1)\int_{s_1}^\infty h(t,s_2) ...\int_{s_{n-1}}^\infty h(t,s_n)\,ds_n...ds_1\right.$$
$$\left.-\int_x^\infty \tilde h(s_1)\int_{s_1}^\infty \tilde h(s_2) ...\int_{s_{n-1}}^\infty \tilde h(s_n)\,ds_n...ds_1\right|
$$
$$\leq C(u(1))^n\left(1+\frac t{x^2}\right)^{n-1}\left(1+\frac 1x\right)\left(\frac{\sqrt{t}}{x}+t^{\frac 16^-}\right),$$
So we define $a_j(x)$ to be $a_j(t,x)$ with the functions $h(t,s)$ replaced by $\tilde h(s)$ and we get for $x\neq 0$:
$$T(t,x)-\sum_{j=1}^{\infty}a_j(x)=\mathcal O\left(t^{\frac 16^-}\right).$$
In particular $T(t,x)$ has a limit at $t=0$, 
$$T(0,x)=\sum_{j=1}^{\infty}a_j(x).$$ 
 The integral equation that $T(0,x)$ satisfies allows us to get a a self-similar estimate on $t\leq x^2$ (see Proposition 3.6 in \cite{BV4}):
$$|T(t,x)-T(0,x)|\leq C(u(1))\,\frac{\sqrt{t}}{x}.$$

\subsection{IVP and continuation after singularity time}\label{sect:IVP}
%%%%%%%%%%%%%%%%%%%%%%%%%%%%%%%%%%%%%%%%%%%%%%
We start by finding some extra informations about the singularity formation for $\chi$ at time $t=0$, that will also indicate us how to construct solutions with initial data curves with a corner. Eventually we shall extend the solution $\chi$ to negative times.
\subsubsection{Formation of the angle}
The full details concerning this subsection can be found in \S5 in \cite{BV3}. 
We denote
$$T_n(x)=T(t_n,\sqrt{t_n}x)\quad,\quad N_n(x)=N(t_n,\sqrt{t_n}x),$$
for a sequence $t_n$ of times that tend to zero. From \eqref{derKoiso} we have
$$T_n'(x)=\sqrt{t_n}\,\Re\,(\overline{\psi}(t_n,\sqrt{t_n}x)N_n(x))=\Re\,( ae^{i\frac{x^2}{4}}N_n(x))+o(t_n)N_n(x),$$
$$N_n'(x)=-\sqrt{t_n}\,\psi(t_n,\sqrt{t_n}x)T_n(x)=-ae^{i\frac{x^2}{4}}T_n(x)+o(t_n)T_n(x).$$
It follows that $\mathcal A=\{T_n,n\in\mathbb N\}$ is a collection of pointwise bounded and equicontinuous functions. Then Arzela-Ascoli theorem allows us to obtain a subsequence, re-called $T_n$, that converges uniformly on any compact subset of $\mathbb R$. We can do the same for $\mathcal B=\{N_n,n\in\mathbb N\}$ and conclude that 
$$\lim_{n\rightarrow\infty}(T_n(x),N_n(x))=(T_o(x),N_o(x)).$$
and the system satisfied by $(T_o(x),N_o(x))$ is 
$$\left\{\begin{array}{c}T_o'(x)=\Re \,(a e^{i\frac{x^2}{4}}N_o(x)),\\ N_o'(x)=ae^{i\frac{x^2}{4}}T_o(x).
\end{array}\right.$$
 Therefore $(T_o(x),\Re \,(e^{-i\frac{x^2}{4}}N_o(x)),\Im\,(e^{-i\frac{x^2}{4}}N_o(x)))$ is the Frenet frame of a curve with curvature and torsion $(a,\frac x 2)$, 
exactly as the profile of the self-similar solution $\chi_a$.

 Hence on the one hand, modulo a rotation $R$, from \cite{GRV}) we have
$$T_o(x)=A_a^++\mathcal O\left(\frac{1}{x}\right)\quad,\quad N_o(x)=B_a^++\mathcal O\left(\frac{1}{x}\right).$$
On the other hand using the self-similar convergence of $T(t)$ to $T(0)$,
$$T_o(x)=\lim_{n\rightarrow\infty}T_n(x)$$
$$=\lim_{n\rightarrow\infty}(T(t_n,\sqrt{t_n}x)-T(0,\sqrt{t_n}x)+T(0,\sqrt{t_n}x))= \mathcal O\left(\frac{1}{x}\right)+\lim_{n\rightarrow\infty}T(0,\sqrt{t_n}x).$$
Therefore we obtain the existence and the value of $T(0,0^\pm)$, modulo the rotation $R$:
$$T(0,0^\pm)=A_a^\pm,$$
so a corner with the same angle as the one of the self-similar solution $\chi_a$ is generated at $t=0, x=0$.
\subsubsection {The initial value problem}
In a similar way we prove that, modulo the rotation $R$, there exists a limit $\tilde N(0,x)$ of $\tilde N(t,x)$ and
$\tilde N(0,0^\pm)=B_a^\pm.$
Moreover,
$$\left\{\begin{array}{c}
T_x(0,x)=-\Re (\frac{\widehat{u_+}\left(\frac x2\right) }{|x|^{ia^2}}\tilde N(0,x)),\,\tilde N_{x}(0,x)=-\overline{\frac{\widehat{u_+}\left(\frac x2\right) }{|x|^{ia^2}}}T(0,x),\\(T,\tilde N)(0,0^\pm)=(A^\pm,B^\pm).\end{array}\right.$$\\
In particular, the recipe for constructing a solution $\chi(t,x)$ with initial data $\chi_0(x)$ displaying a corner as $\chi_a(0,x)$ and curvature in weighted spaces as stated in Theorem \ref{theivp} is the following:
\begin{itemize}
\item we define $g(x)$ and $\tilde N_0(x)\perp T_0(x)$ for $x>0$ (similarly for $x<0$) by
$$
T_{0x}(x)=-\Re (g(x)\tilde N_0(x)),\,\tilde N_{0x}(x)=-\overline{g}(x)T_0(x),\, (T_0,\tilde N_0)(0)=(A^+,B^+),
$$
\item $|g(x)|$ is the curvature of $\chi_0$, and we define the $X^\gamma$-type function
$$u_+(x)=\mathcal F^{-1}\left(g(2\cdot)|2\cdot|^{ia^2}\right),$$
\item from the existence of wave operator with final data $u_+$ we obtain a solution $u(t,x)$ of \eqref{NLSu},
\item from $u(t,x)$ we obtain, as explained in the first section, a solution $\chi$ of the binormal flow.
\end{itemize}
Now, modulo a rotation, $\chi$ has initial data $\chi_0$. Indeed, by doing the analysis in \S \ref{subsect:tang} and above in the current subsection \S \ref{sect:IVP}, it follows that the tangent of $\chi(t,x)$ has a trace $(T,\tilde N)(0,x)$ at $t=0$, that satisfies, modulo a rotation $R$, the same system as $(T_0,\tilde N_0)$, so $T_0=RT(0)$. Then, one shows that indeed $T(0)=\partial_x\chi(0)$. By translating $R\chi$ by $\chi_0(0)-R\chi(0,0)$ we obtain a solution of \eqref{VFE} with initial data $\chi_0$. In \S 3.5 in \cite{BV4} we show by a Gronwall argument that the uniqueness result holds. 

%For the self-similar solutions, a way to continue them after $t=0$ is by the only unique binormal flow solution with initial data $\chi_a(0,-x)$, given by $\pi-$rotation with respect to the bisector of the angle. remark on time reversibility, continuation after singularity time for the self-similar solutions and their $A^\pm$,$B^\pm$ vectors, 

\subsubsection{Continuation after singularity time}
From the reversibility of \eqref{VFE} ($\chi(t,x)$ solution implies $\chi(-t,-x)$ solution also), continuing the solution $\chi(t,x)$ for negative times means solving the initial value problem for positive times with initial data $\chi_0^\star(x)=\chi(0,-x)$, oriented curve with the $\chi_a(0,-x)$ corner. Now, the unique self-similar solution with initial data $\chi_a(0,-x)$ is the $\pi-$rotation $R\chi_a(t,x)$ around the bisector of the angle. Then the corresponding asymptotic vectors are $R A^\pm,RB^\pm$, that can be written as $R^\mp(-A^\mp,\overline{B^\mp})$ for some particular rotations $R^\mp$ (see Proposition 2.2 in \cite{BV4}).

So the IVP recipe starts with initial data $(T_0^*,\tilde N_0^*)(0)=\tilde R^\mp(-A^\mp,\overline{B^\mp})$ and in particular we have as solutions
$$N_0^\star(x)=\tilde R^\mp\overline{\tilde N(0,-x)},\quad g^\star(x)=\overline{g(-x)}=\widehat{\overline{u_+}}\left(\frac x2\right) |x|^{ia^2}.$$
Therefeore the IVP result gives us a solution $\chi^\star(t,x)$ for positive times with initial data $\chi_0^\star$. Then $\chi(t,x)=\chi^\star(-t,-x)$ is the extension for negative times of $\chi(t,x)$.

\bibliographystyle{plain}

\begin{thebibliography}{123}

\bibitem{ArHa} R.J. Arms and F.R. Hama, 
{\em Localized-induction concept on a curved vortex and motion of an elliptic vortex
ring}, Phys. Fluids, (1965), 553.




\bibitem{BV1}  V.~Banica and L.~Vega, 
{\it On the stability of a singular vortex dynamics}, 
Comm. Math. Phys. {\bf 286} (2009), 593--627.

\bibitem{BV2}  V.~Banica and L.~Vega, 
{\it Scattering for 1D cubic NLS and singular vortex dynamics}, 
J. Eur. Math. Soc. {\bf 14} (2012), 209--253.

\bibitem{BV3}  V.~Banica and L.~Vega, 
{\it Stability of the self-similar dynamics of a vortex filament}, 
to appear in Arch. Ration. Mech. Anal.

\bibitem{BV4}  V.~Banica and L.~Vega, 
{\it The initial value problem for the binormal flow with rough data}, 
ArXiv:1304.0996. 

\bibitem{Ben-Cag-Mar} D.~Benedetto, E.~Caglioti, and C.~Marchioro, 
 {\it On the motion of a vortex ring  with a sharply concentrated vorticity},  
Math. Methods Appl. Sci. {\bf 23} (2000), 147--168.

\bibitem{BeOrSm} 
F.~B\'ethuel, G.~Orlandi, and D.~Smets, 
 {\it Convergence of the parabolic Ginzburg-Landau equation to motion by mean curvature}, 
Ann. of Math. {\bf 163} (2006), 37--163.

\bibitem{Bo} J~ Bourgain,
{\it Fourier transform restriction phenomena for certain lattice subsets and applications to nonlinear evolution equations I. Schr�dinger equations,} 
Geom. Funct. Anal. {\bf  3} (1993), 107--156.

\bibitem{Bu} T. F. Buttke, {\it A numerical study of superfluid turbulence in the Self Induction Approximation}, 
J. of Compt. Physics {\bf 76} (1988), 301--326


\bibitem{Cal} A. Calini and T. Ivey, {\it Stability of Small-amplitude Torus Knot Solutions of the Localized Induction Approximation}, 
J. Phys. A: Math. Theor. 44 (2011) 335204.


\bibitem{Ca} 
R.~Carles, 
{\it Geometric Optics and Long Range Scattering for One-Dimensional Nonlinear Schr\"odinger Equations, }
Comm. Math. Phys. 
{\bf 220} 
(2001),  
41--67.


\bibitem{CW} 
T.~Cazenave and F.B.~Weissler,
{\it The Cauchy problem for the critical nonlinear Schr�dinger equation,} 
Non. Anal. TMA {\bf 14} (1990), 807--836. 

\bibitem{CCT}
M.~Christ, J.~Colliander, and T.~Tao, 
{\it Ill-posedness for nonlinear
  {S}chr\"odinger and wave equations}, 
ArXiv:0311048.



\bibitem{DaR} L. S. Da Rios,  
{\em On the motion of an unbounded fluid with a vortex filament of any shape},  
Rend. Circ. Mat. Palermo {\bf 22} (1906), 117.

\bibitem{GMS} 
P.~Germain, N.~Masmoudi and J.~Shatah, 
{\it Global solutions for 2D quadratic Schr�dinger equations,} 
J. Math. Pures Appl. {\bf 97}(2012), 505--543.

\bibitem{GV}
J.~Ginibre and G.~Velo, 
{\it On a class of nonlinear {S}chr\"odinger
  equations. {II} {S}cattering theory, general case}, 
  J. Funct. Anal. {\bf 32} (1979), 33--71.
  
  \bibitem{Gr}
A.~Gr\"unrock, 
{\it Bi- and trilinear Schr?dinger estimates in one space dimension with applications to cubic NLS and DNLS, }
  Int. Math. Res. Not.  
  {\bf 41}
 (2005), 2525--2558.


\bibitem{GNT}  S. Gustafson, K. Nakanishi, T.-P. Tsai, 
Global dispersive solutions for the Gross-Pitaevskii equation in two and three dimensions, 
 Ann. Henri Poincar\'e  8  (2007),  no. 7, 1303--1331. 


\bibitem{GRV} S.~Guti\'errez, J.~Rivas and L.~Vega,  
{\it Formation of singularities and self-similar vortex motion under the localized induction approximation}, 
Comm. Part. Diff. Eq. {\bf 28} (2003), 927--968.

\bibitem{Ha}
H.~Hasimoto,
 {\it A soliton in a vortex filament, }
J. Fluid Mech.
{\bf 51}
(1972), 477--485.

\bibitem{HN} 
N.~Hayashi and P.~Naumkin, 
{\it Domain and range of the modified wave operator for Schr?odinger equations
with critical nonlinearity,} 
Comm. Math. Phys. {\bf 267} (2006), 477--492.

\bibitem{HoBr}
E.J. Hopfinger, F.K. Browand,
{\it Vortex solitary waves in a rotating, turbulent flow, }
Nature
{\bf 295}, 
(1981), 
393--395.

\bibitem{HGCV}
F.~de~la~Hoz, C.~Garc\'{i}a-Cervera and L.~Vega,
{\it A numerical study of the self-similar solutions of the Schr\"odinger Map}, 
SIAM J. Appl. Math. {\bf 70} (2009), 1047--1077. 

\bibitem{Je} R.~L.~Jerrard, 
{\it Vortex filament dynamics for Gross-Pitaevsky type equations}, Ann. Scuola Norm. Sup. Pisa Cl. Sci. {\bf 5} (2002), 733--768.

\bibitem{JeSm1}
R.~L.~Jerrard and D. Smets,
{\em On Schr\"odinger maps from $T^1$ to $S^2$,}
arXiv:1105.2736.

\bibitem{JeSm2}
R.~L.~Jerrard and D.~Smets,
{\em On the motion of a curve by its binormal curvature,}
arXiv:1109.5483.

\bibitem{KPV}
C. Kenig, G. Ponce and L. Vega,
 {\em On the ill-posedness of some canonical non-linear dispersive equations, }
Duke Math. J.
{\bf 106} 
(2001)
716--633.

\bibitem{Ko} N.~Koiso,
 {\em Vortex filament equation and semilinear Schr\"odinger equation, }
Nonlinear Waves, Hokkaido University Technical Report Series in Mathematics
{\bf 43} (1996)
221--226.


\bibitem{Laf} S. Lafortune, {\em Stability of solitons on vortex filaments,}
Phys. Lett. A {bf 377} (2013), 766--769.

\bibitem{LD} M.~Lakshmanan and M.~Daniel,
{\em On the evolution of higher dimensional Heisenberg continuum spin systems}, Physica A (1981), 107, 533--552.

\bibitem{LRT} M.~Lakshmanan, T.~W.~Ruijgrok and C.~J.~Thompson,
{\em On the the dynamics of a  continuum spin system}, Physica A (1976), 84, 577--590.

\bibitem{LC} T.~Levi-Civita, 
{\it Attrazione Newtoniana dei Tubi Sottili e Vortici Filiformi}, 
Ann. Scuola Norm. Sup. Pisa Cl. Sci. {\bf 1} (1932), 229--250

\bibitem{Li} 
T.~Lipniacki, 
{\it Quasi-static solutions for quantum vortex motion under the localized induction approximation,}
J. Fluid Mech. {\bf 477} (2002), 321--337.

\bibitem{Mag}
F. Maggioni, S. Z. Alamri, C. F. Barenghi, and R. L. Ricca, {\it Velocity, energy and helicity of vortex knots and unknots},  
Phys. Rev. E {\bf 82} (2010), 26309--26317.


\bibitem{MaBe} 
A.~Majda and A.~ Bertozzi, 
{\it Vorticity and incompressible flow.,} 
Cambridge Texts in Applied Mathematics, 27. Cambridge University Press, Cambridge, 2002. 


\bibitem{Mar-rings}
 C.~Marchioro, {\it Large smoke rings with concentrated vorticity}, 
J. Math. Phys. {\bf 40} (1999), 869--883.


\bibitem{Mar-Neg} C.~Marchioro and P.~Negrini, 
{\it On a dynamical system related to fluid
mechanics}, 
NoDEA Nonlinear Differential Equations Appl. {\bf 6}  (1999), 473--499.

\bibitem{MTT} 
K.~Moriyama, S.~Tonegawa and Y.~ Tsutsumi, 
{\it Wave operators for the nonlinear Schr�dinger equation with a nonlinearity of low degree in one or two space dimensions,}
Commun. Contemp. Math. {\bf 5} (2003), 983--996 .

\bibitem{NiTa}
T.~Nishiyama and A.~Tani, 
{\it Solvability of the localized induction equation for vortex motion,} 
Comm. Math. Phys. {\bf 162} (1994), 433?-445.

\bibitem{Oz} 
T.~Ozawa, 
{\it Long range scattering for nonlinear Schr\"odinger equations in one space dimension, }
Comm. Math. Phys. {\bf 139} 
(1991), 
479--493. 

\bibitem{PMQ}
C.~S.~Peskin and D.~M.~McQueen, 
{\it Mechanical equilibrium determines the fractal fiber architecture of aortic heart valve leaflets}, 
Am. J. Physiol. 266 (Heart Circ.
Physiol. 35) (1994), H319--H328.


\bibitem{Ri1}  R.~L. Ricca, 
 {\em The contributions of Da Rios and Levi-Civita to asymptotic potential theory and vortex filament dynamics}, 
Fluid Dynam. Res. {\bf 18} (1996), 245--268. 

\bibitem{Ri2} R.L.~Ricca, 
{\it Rediscovery of Da Rios equations,}
Nature {\bf 352} (1991), 561--562.


\bibitem{Sc}
K.W.~Schwarz, 
{\it Three-dimensional vortex dynamics in superfluid $^4$He: Line-line and line-boundary
interactions}, Phys. Rev B {\bf 31} (1985), 5782--5804.

\bibitem{ST}
A.~Shimomura and S.~Tonegawa, 
{\it Long-range scattering for nonlinear Schr�dinger equations in one and two space dimensions,} 
Differ. Integral Equ. {\bf 17} (2004), 127--150.


\bibitem{TN} 
A.~Tani and T.~Nishiyama, 
{\it Solvability of equations for motion of a vortex filament with or without axial flow,} 
Publ. Res. Inst. Math. Sci. {\bf 33} (1997), 509?-526.

\bibitem{VV}
A.~Vargas and L.~Vega,
 {\it Global well-posedness for 1d non-linear  Schrodinger equation for data with an infinite $L^2$ norm, }
J. Math. Pures Appl.
{\bf 80}  (2001), 
1029--1044.

\bibitem{V}
E.J.~Vigmond, C.~Clements, D.M.~McQueen and C.S.~Peskin, 
{\it Effect of bundle branch block on cardiac output: A whole heart simulation study,} 
Prog. Biophys. Mol. Biol. {\bf 97} (2008), 520--42.

\end{thebibliography}

\end{document}